\documentclass[aos,MSNbibl,seceqn,nameyear,dvips]{arximspdf}

\doi{10.1214/11-AOS930}
\volume{39}
\issue{6}
\pubyear{2011}
\firstpage{3003}
\lastpage{3031}

\makeatletter

\newcommand{\iint}{\int\!\!\int}

\newtheorem{theorem}{Theorem}[section]
\newtheorem{lem}{Lemma}[section]
\newtheorem{cor}{Corollary}[section]

\newproclaim{assum}{Assumption}[section]
\newproclaim{exm}{Example}[section]
\newproclaim{remark}{Remark}[section]

\makeatother

\begin{document}
\begin{frontmatter}

\title{Posterior consistency of nonparametric conditional moment
restricted models}
\runtitle{Conditional moment restricted models}

\begin{aug}
\author[A]{\fnms{Yuan} \snm{Liao}\corref{}\ead[label=e1]{yuanliao@princeton.edu}}
\and
\author[B]{\fnms{Wenxin} \snm{Jiang}\ead[label=e2]{wjiang@northwestern.edu}}
\runauthor{Y. Liao and W. Jiang}
\affiliation{Princeton University and Northwestern University}
\address[A]{Department of Operations Research\\
\quad and Financial Engineering\\
Princeton University\\
Sherrerd Hall\\
Princeton, New Jersey 08544\\
USA\\
\printead{e1}}
\address[B]{Department of Statistics\\
Northwestern University\\
2006 Sheridan Rd\\
Evanston, Illinois 60208\\
USA\\
\printead{e2}} 
\end{aug}

\received{\smonth{3} \syear{2010}}
\revised{\smonth{9} \syear{2011}}

%
\begin{abstract}
This paper addresses the estimation of the nonparametric conditional
moment restricted model that involves an infinite-dimensional parameter
$g_0$. We estimate it in a \textit{quasi-Bayesian} way, based on the
limited information likelihood, and investigate the impact of three
types of priors on the posterior consistency: (i) truncated prior
(priors supported on a bounded set), (ii) thin-tail prior (a prior that
has very thin tail outside a growing bounded set) and (iii) normal
prior with nonshrinking variance. In addition, $g_0$ is allowed to be
only partially identified in the frequentist sense, and the parameter
space does not need to be compact. The posterior is regularized using a
slowly growing sieve dimension, and it is shown that the posterior
converges to any small neighborhood of the identified region. We then
apply our results to the nonparametric instrumental regression model.
Finally, the posterior consistency using a random sieve dimension
parameter is studied.
\end{abstract}

%

\end{frontmatter}

\section{Introduction}\label{sec1}
We consider a conditional moment restricted model
%
%
\begin{equation}\label{equ11}
E(\rho(Z, g_0)|W, g_0)=0,
\end{equation}
where $(Z^T,W^T)$ is a vector of observable random variables, and $W$
may or may not be included in $Z$. Here $\rho$ is a one-dimensional
residual function known up to~$g_0$. The conditional expectation is
taken with respect to the conditional distribution of $Z$ given $W$ and
$g_0$, assumed unknown. The parameter of interest is~$g_0$, which is
infinite dimensional. Moreover, suppose we observe independent and
identically distributed data $\{(Z_i^T, W_i^T)\}_{i=1}^n$ of~$(Z^T,W^T)$.

Model (\ref{equ11}) is a very general setting, which encompasses many important
classes of nonparametric and semiparametric models.
%
%
\begin{exm}[(Regular nonparametric regression)]\label{examp11}
Consider the model
\[
Y=g_0(W)+\varepsilon
\]
assuming $E(\varepsilon|W)=0$. Let $Z=(Y,W)$, then it can be written as
the conditional moment restricted model with $\rho(Z,g_0)=Y-g_0(W)$.
\end{exm}
%
%
\begin{exm}[(Single index model)]\label{examp12}
Consider the single index model
\[
Y=h_0(W^T\theta_0)+\varepsilon,
\]
where $E(\varepsilon|W)=0$. The parameter of interest is $(h_0, \theta
_0)$, with $h_0$ being nonparametric. This type of model is studied by
\citet{Ich93} and \citet{AntGreMcK04}. By defining $Z=(Y, W)$,
$g_0=(h_0, \theta_0)$ and $\rho(Z,g_0)=Y-h_0(W^T\theta_0)$, we can
write $E(\rho(Z,g_0)|W, g_0)=0$.
\end{exm}
%
%
\begin{exm}[(Nonparametric IV regression)]\label{examp13}
Consider the nonparametric model
\[
Y=g_0(X)+\varepsilon,
\]
where $X$ is an endogenous regressor, meaning that $E(\varepsilon|X)$ does
not vanish. However, suppose we have observed an instrumental variable
$W$ for which $E(\varepsilon|W)=0$; then it becomes a nonparametric
regression model with instrumental variables (NPIV), studied by
\citet{NewPow03} and \citet{HalHor05}. Define $\rho
(Z,g_0)=Y-g_0(X)$, with $Z=(Y,X)$. Then we have the conditional moment
restriction.
\end{exm}
%
%
\begin{exm}[(Nonparametric quantile IV regression)]\label{examp14} The
nonparametric
quantile IV regression was previously studied by \citet{CheHan05},
\citet{CheImbNew07} and \citet{HorLee07}. The
model is
\[
y=g_0(X)+\varepsilon,\qquad P(\varepsilon\leq0|W)=\gamma,
\]
where $g_0$ is the unknown function of interest, and $\gamma\in(0,1)$
is known and fixed. Assume $X$ is a continuous random variable. Then
the conditional moment restriction is given by
\[
E(\rho(Z,g_0)|W, g_0)=0,\qquad \rho(Z,g_0)=I_{(y\leq
g_0(X))}-\gamma.
\]
\end{exm}

If we define $G(g)=E_W[E(\rho(Z,g)|W, g_0)]^2$, an equivalent way of
writing model (\ref{equ11}) is then $G(g_0)=0$. When the unknown
function $g_0$
depends on certain endogenous variable as in Examples \ref{examp13} and
\ref{examp14}, the
identification and consistent estimation of $g_0$ is challenging. On
one hand, there can be multiple functions in the parameter space that
satisfy the moment restriction~(\ref{equ11}). On the other hand, even
if $g_0$
is identified, [in which case the functional $G(g)$ is uniquely
minimized at $g=g_0$, as is typically assumed in the literature],
reducing $G(g)$ toward $G(g_0)$ does not guarantee that $\|g-g_0\|_s$
will also be close to zero, for a certain norm \mbox{$\|\cdot\|_s$} of
interest. %
Therefore, minimizing a consistent estimator of $G(g)$ does not lead to
a consistent estimator of $g_0$ under \mbox{$\|\cdot\|_s$}. This phenomenon is
usually known as the ``ill-posed inverse problem'' in the literature.

The general form of (\ref{equ11}) was first studied by \citet
{AiChe03} and
\citet{NewPow03}, where the authors considered sieve
approximation of $g_0$ and estimated it in a compact parameter space.
Recently, \citet{ChePou} relaxed the compactness assumption and
achieved the consistency and convergence rate using the penalized sieve
minimum distance estimation. In recent years there has also been
extensive literature on the NPIV model (Example \ref{examp13}) itself.
In these
papers, the authors introduce a~Tikhonov tuning parameter to play a
role of ``regularization'' in order to overcome the ill-posed inverse
problem; see, for example, \citet{HalHor05} and \citet{Daretal}.
Other related works on the nonparametric
instrumental variables can be found in \citet{CheGagSca},
\citet{JohVanVan},
Horowitz (\citeyear{Hor07}, \citeyear{Hor11}), among others.

Compared to the growing literature from the frequentist perspective,
there is very little understanding of the consistent estimation using
either a Bayesian or a quasi-Bayesian approach. This paper proposes a
quasi-Bayesian procedure and studies the impact of various priors of
$g_0$ on the posterior consistency. Our setup is built on a sieve
approximation technique similar to \citet{ChePou}, which
assumes that $g_0$ can be approximated arbitrarily well on a
finite-dimensional sieve space. In order to keep our procedure robust
to the
distribution specification and convenient for practical implementation,
without specifying a known distribution on the data generating process,
we employ a limited information likelihood [\citet{Kim02} and
\citet{LiaJia10}], a moment-condition-based Gaussian approximated
likelihood. The use of such a likelihood is more straightforward for
models characterized by either moment conditions or estimating
equations than the common methods based on Dirichlet process priors in
the nonparametric Bayesian literature. With priors placed directly on
the sieve coefficients, we show that the proposed posterior is
consistent. Due to the difficulty of identifying $g_0$ in practice, we
do not assume $g_0$ to be necessarily identified. As a result the
posterior consistency here means that, asymptotically, the posterior
converges into arbitrarily small neighborhood of the region where $g_0$
is partially identified. Therefore, we also extend model (\ref{equ11})
to the
partial identification setup [\citet{CheHonTam07}
and \citet{San}]. We will consider three types of priors: (i) priors
supported on a bounded set (truncated prior), (ii) priors with tails
decaying fast outside a bounded set (thin-tail prior) and (iii)
Gaussian priors with nonshrinking variance.

Recently, \citet{FloSimN1} proposed a quasi-Bayesian
approach for the NPIV model. They assumed that the error term follows a
normal distribution and achieved consistency by regularizing an
operator that defines the posterior mean. Our approach differs from
theirs essentially in the way of overcoming the ill-posed inverse
problem. While \citet{FloSimN1} put a Gaussian prior on an
infinite-dimensional function space, they require the variance of the
prior to shrink to zero. In contrast, we place the prior directly on
the sieve coefficients in a finite-dimensional vector space and require
the sieve dimension to grow slowly with the sample size. Our approach
then corresponds to Chen and Pouzo's (\citeyear{ChePou}) sieve minimum distance
procedure using slowly growing sieves. As a result, it is the
finite-dimensional sieve that plays the role of regularization instead
of a
shrinking prior. In addition, our approach allows nonnormal priors.

Models based on moment conditions as (\ref{equ11}) have been proved to be
essential in many statistical applications, such as financial asset
pricing [\citet{GalTau89}, \citet{CheLud09}],
consumer behavior in economics [\citet{BluCheKri07},
\citet{San}] and return to college education [\citet{Hor11}].
Therefore, this paper develops a quite convenient and straightforward
quasi-Bayesian approach for these applied problems.

The remainder of this paper is organized as follows: Section \ref{sec2}
introduces general theorems on two types of posterior consistency,
which provide sufficient conditions under which a posterior constructed
on a sieve space is consistent. Section~\ref{sec3} specifies the priors and
shows the consistency results by verifying the sufficient conditions
given in Section \ref{sec2}. Section \ref{sec4} studies in detail the
NPIV model as a
specific example. Section \ref{sec5} discusses the case of the random sieve
dimension. Finally, Section \ref{sec6} concludes with further discussions.
Proofs are given in the supplementary material.

Throughout the paper, for any two positive deterministic sequences $\{
a_n\}_{n=1}^{\infty}$ and $\{b_n\}_{n=1}^{\infty}$, write $a_n\succ
b_n$ and $b_n\prec a_n$ if $b_n=o(a_n)$.
In addition, $a_n\sim b_n$ if there exist
$c_1$ and $c_2>0$ such that $c_1b_n\leq a_n\leq c_2b_n$
for all large enough $n$.

\section{General posterior consistency theorems}\label{sec2}
\subsection{Sieve approximation}
Suppose we are interested in a nonparametric regression function
$g_0\in
(\mathcal{H},\mbox{$\|\cdot\|_s$})$. which is assumed to be inside an
infinite-dimensional Banach space $\mathcal{H}$ endowed with norm $\|
\cdot\|_s$.
Examples of the space $(\mathcal{H}, \mbox{$\|\cdot\|_s$})$ include: space of
bounded continuous functions with norm $\|g\|_s=\sup_x|g(x)|$, the
space of square integrable functions $\{g\dvtx E[g(X)^2]<\infty\}$ with
$\|
g\|_s=\sqrt{E[g(X)^2]}$, etc. In addition, suppose there exists a set
of basis functions $\{\phi_1,\phi_2,\ldots\}\subset\mathcal{H}$ such that
$g_0\in\mathcal{H}$ can be approximated by a truncated sum $g_b=\sum
_{i=1}^{q_n}b_i\phi_i$ for a vector of coefficients
$(b_1,\ldots,b_{q_n})^T$, where $q_n$ is a pre-determined constant that
grows to infinity. Then $g_b$ lies in an approximating space $\mathcal
{H}_n$ spanned by $\{\phi_1,\ldots,\phi_{q_n}\}$. Here $\mathcal{H}_n$
grows to be dense in $\mathcal{H}$, called a sieve approximating space.

There is extensive literature on the posterior consistency using
sieve approximation. \citet{SheWas01} applied an orthogonal
basis expansion to the nonparametric regression problem. \citet{Wal03}
and \citet{ChoSch07} provided general results for a class of
Bayesian regression models when the data have a normal distribution.
Other results on nonparametric regression problems can be found, for
example, in \citet{Hua04}, \citet{Ghovan07}, etc.

Suppose we are given $n$ independent identically distributed
observations $X^n=(X_1,X_2,\ldots,X_n)$. In this paper we do not assume
any specific distribution of $X^n|g_0$, but propose a quasi-Bayesian
approach, which is based on a pseudo-likelihood,
\[
L(g_b)=\exp\biggl(-\frac{n}{2}\bar{G}(g_b)\biggr),
\]
where $\bar{G}\dvtx\mathcal{H}_n\rightarrow[0,\infty)$ is a stochastic
functional, which we call the \textit{sample risk functional}. Suppose
there exists a nonnegative functional $G$, such that for a bounded set
$\mathcal{F}_n\subset\mathcal{H}_n$,
\[
{\sup_{g_b\in\mathcal{F}_n}}|\bar{G}(g_b)-G(g_b)|=o_p(1).
\]
We call $G$ the \textit{objective functional} or \textit{risk
functional} throughout the paper.

In the literature, it is often assumed that the true regression
function $g_0$ is point identified (as opposed to ``partially
identified'' in the following) as the unique minimizer of $G$ on
$\mathcal{H}$, that is,
\[
\{g_0\}=\mathop{\arg\min}_{g\in\mathcal{H}}G(g).
\]
Then quasi-Bayesian approaches usually construct $\bar{G}$ as the
sample analog of~$G$ [see \citet{CheHon03}]. In many applications of the model considered in
this paper, however, it is more natural to assume that $G$ has multiple
global minimizers on $\mathcal{H}$; see detailed discussions in Section
\ref{sec3}. In this case, we say $g_0$ is \textit{partially identified}
(in the
frequentist sense) on
\[
\Theta_I=\mathop{\arg\min}_{g\in\mathcal{H}}G(g),
\]
and $\Theta_I$ is called the \textit{identified region}. Therefore
$\Theta_I$ is the main object of interest in this paper.

For any $b=(b_1,\ldots,b_{q_n})^T\in\mathbb{R}^{q_n}$, let $g_b=\sum
_{i=1}^{q_n}b_i\phi_i$. Similarly to the standard treatments in
\citet{SmiKoh96} and
\citet{AntGreMcK04},\vadjust{\goodbreak} we put prior $\pi(b)$ on
the sieve coefficients $b=(b_1,b_2,\ldots,\allowbreak b_{q_n})$, and obtain a
posterior distribution,
\[
P(g_b|X^n)\propto\pi(b)L(g_b).
\]
For any $g_1\in\mathcal{H}$, define
\[
d(g_1,\Theta_I)={\inf_{g\in\Theta_I}}\|g_1-g\|_s,
\]
and the \textit{$\varepsilon$-expansion} as a neighborhood of the
identified region
\[
\Theta_I^{\varepsilon}=\{g\in\mathcal{H}\dvtx d(g,\Theta_I)<\varepsilon
\}.
\]
Then the posterior consistency in this paper refers to the following:
for any $\varepsilon>0$,
\[
P(g\in\Theta_I^{\varepsilon}|X^n)\rightarrow^p1.
\]

\subsection{Posterior consistency theorems}
We first present two theorems of general posterior consistency using
the sieve approximation, which involve conditions on the tail
probability of $\pi$ as well as the performance of $\bar{G}$. They are
based on the following variant of an inequality from \citet{JiaTan08},
Proposition 6. These inequalities will be proved in the
supplementary material [\citet{LiaJiaN1}]:
%
%
\begin{lem}\label{l21}
Suppose the support of the prior $\pi$ can be partitioned as $\mathcal
{F}_n\cup\mathcal{F}_n^c$. Then for any deterministic sequence
$\delta_n>0$,
%
%
\begin{eqnarray}\label{e21}
&&E\Bigl\{P\Bigl(G(g_b)-\inf_{g\in\mathcal{H}}G(g)>5\delta_n|X^n\Bigr)\Bigr\}\nonumber\\
&&\qquad\leq P\Bigl({\sup
_{g\in\mathcal{F}_n}}|\bar{G}(g)-G(g)|\geq\delta_n\Bigr)\nonumber\\[-8pt]\\[-8pt]
&&\qquad\quad{}
+\frac{e^{-2n\delta_n}}{\pi(G(g_b)-\inf_{g\in\mathcal
{H}}G(g)<\delta_n
\cap g_b\in\mathcal{F}_n)}\nonumber\\
&&\qquad\quad{}+EP(g_b\in\mathcal{F}_n^c|X^n).\nonumber
\end{eqnarray}
In addition,
\begin{eqnarray*}
EP(g_b\in\mathcal{F}_n^c|X^n)&\leq& P\Bigl({\sup_{g\in\mathcal
{F}_n}}|\bar
{G}(g)-G(g)|\geq\delta_n\Bigr)\\
&&{}+\frac{\pi(\mathcal{F}_n^c)e^{2n\delta_n}}{\pi(G(g_b)-\inf
_{g\in
\mathcal{H}}G(g)<\delta_n \cap g_b\in\mathcal{F}_n)}.
\end{eqnarray*}
\end{lem}

These inequalities imply the following result on the risk consistency:
%
%
\begin{theorem}[(Risk consistency)] \label{t21} Suppose the following
conditions hold with respect to a deterministic positive sequence
$\delta_n$:

\begin{longlist}
\item
Tail condition: as $q_n$ and $n\rightarrow\infty$, either
$EP(g_b\in
\mathcal{F}_n^c|X^n)=o(1)$ or
$\pi(\mathcal{F}_n^c)=O(e^{-4n\delta_n})$.\vspace*{1pt}

\item Approximation condition: $\pi(G(g_b)-\inf_{g\in\mathcal
H}G(g)<\delta_n, g_b\in\mathcal{F}_n)\succ e^{-2n\delta_n}$.

\item Uniform convergence: $P[{\sup_{g\in\mathcal{F}_n}}|\bar
G(g)-G(g)|\geq
\delta_n]=o(1)$.\vspace*{1pt}

Then we have the risk consistency result at rate $\delta_n$
\[
P\Bigl(G(g_b)-\inf_{g\in\mathcal{H}} G(g)<\delta_n|X^n\Bigr)=1-o_p(1).
\]
\end{longlist}
\end{theorem}

The naming of these conditions is obvious, except for (ii). There, the
approximation refers to the ability of the functions in $\mathcal{ F}_n$
(proposed by the prior $\pi$) to approximately minimize the risk ${G}$
over $ \mathcal{H}$ with not-too-small prior probability.

When the following condition is added, the risk consistency leads to the
estimation consistency.
%
%
\begin{theorem}[(Estimation consistency)] \label{t22} Suppose there
exists a
sequen\-ce~$\delta_n$ such that the following
conditions hold:\vspace*{8pt}

\textup{(i), (ii), (iii)} in the previous theorem;

\textup{(iv)} (distinguishing ability) for any $\varepsilon>0$,
\[
\inf_{g \in\mathcal{H}_n, g\notin\Theta_I^{\varepsilon}}G(g)-\inf
_{g\in
\mathcal{H}}G(g)\succ\delta_n.
\]
Then for any $\varepsilon>0$, we have
%
%
\begin{equation} \label{e22}
P(g_b \in\Theta_I^{\varepsilon}|X^n)\rightarrow^p1.
\end{equation}
\end{theorem}
\begin{pf}
Theorem \ref{t21} is implied by Lemma \ref{l21}. Now we prove
Theorem~\ref{t22}. For any $\varepsilon>0$, by Theorem \ref{t21},
\begin{eqnarray*}
&&
P(g_b\notin\Theta_I^{\varepsilon}|X^n)\\
&&\quad\leq P\Bigl(g_b\notin\Theta
_I^{\varepsilon
}, G(g_b)-\inf_{g\in\mathcal{H}}G(g)<\delta_n|X^n\Bigr)+o_p(1)\\
&&\quad\leq P\Bigl(g_b\notin\Theta_I^{\varepsilon}, G(g_b)\geq\inf_{g\in
\mathcal
{H}_n, g\notin\Theta_I^{\varepsilon}}G(g), G(g_b)-\inf_{g\in
\mathcal
{H}}G(g)<\delta_n|X^n\Bigr)
+o_p(1)\\
&&\quad\leq P\Bigl(g_b\notin\Theta_I^{\varepsilon}, \delta_n< G(g_b)-\inf
_{g\in
\mathcal{H}}G(g)<\delta_n|X^n\Bigr)+o_p(1)\\
&&\quad=o_p(1),
\end{eqnarray*}
where the third inequality is implied by condition (iv) for all large $n$.
\end{pf}

As a special case of these results, note that when $g_0$ is point
identified as the unique minimizer of $G(g)$ on $\mathcal{H}$, that is,
$\Theta_I=\{g_0\}$, (\ref{e22}) then becomes
\[
P(\|g_b-g_0\|_s<\varepsilon|X^n)\rightarrow^p1,
\]
the regular posterior consistency result.\vadjust{\goodbreak}

In the subsequent sections, we will construct a so-called \textit
{limited information likelihood} $\bar{G}(g)$ and apply the previous
two theorems to the conditional moment restricted model (\ref{equ11}), by
verifying conditions (i)--(iv).

\section{Conditional moment-restricted model}\label{sec3}

\subsection{Limited information likelihood}\label{sec31}
Consider a conditional moment condition
%
%
\begin{equation}\label{e31}
E[\rho(Z, g_0)|W, g_0]=0,
\end{equation}
where $g_0\in\mathcal{H}$ is the true nonparametric structural
function. Here $W$ is $d$-dimensional, with fixed $d$. For simplicity,
throughout the paper, let us assume~$W$ is supported on $[0,1]^d$, as
one can always apply the transformation on each component of $W$,
$W_i\rightarrow\Phi(W_i)$, where $\Phi(\cdot)$ is the standard normal
cumulative distribution function. We focus on the case when $\rho$ is a
one-dimensional function.

Following the setting of \citet{AiChe03} and
\citet{ChePou}, we approximate $\mathcal{H}$ by a sieve space
$\mathcal{H}_n$
that grows to be dense in $\mathcal{H}$. Here $\mathcal{H}_n$ is a
finite-dimensional space spanned by sieve basis functions$\{\phi
_1,\ldots,\phi_{q_n}\}$ such as splines, power series, wavelets and
Fourier series.

As the first step, we transform the conditional moment restriction into
unconditional moment restrictions (but still conditional on $g_0$). Let
$\{[(i-1)/k_n,i/k_n]\}_{i=1}^{k_n}$ be a partition of $[0,1]$, for some
$k_n\in\mathbb{N}$. We then obtain a~partition of the support of $W\dvtx
[0,1]^d=\bigcup_{j=1}^{k_n^d}R_j^n$, where for each $j=1,\ldots,k_n^d$,
%
%
\begin{equation}\label{partition}
R_j^n=\prod_{l=1}^d\biggl[\frac{i_l-1}{k_n}, \frac{i_l}{k_n}\biggr]
\qquad\mbox{for some } i_l\in\{1,\ldots, k_n\}.
\end{equation}
We require $k_n\rightarrow\infty\mbox{ as } n\rightarrow\infty$. Let
$X=(Z,W)$. For each $j$, define
\[
m_{nj}(g, X)=\rho(Z, g)I_{(W\in R_j^n)},
\]
where $I_{(\cdot)}$ is the indicator function. Let $m_n(g, X)=(m_{n1}(g,
X),\ldots, m_{n k_n^d}(g,\allowbreak X))^T$, which is a $k_n^d\times1$ vector.
Equation (\ref{e31}) then implies
%
%
\begin{equation} \label{e33}
Em_n(g_0, X)=0,
\end{equation}
where the expectation is taken with respect to the joint distribution
of $X=(Z,W)$ conditional on $g_0$. Throughout the paper, the
expectation is always taken conditionally on $g_0$. When $k_n>q_n$
there are more moment conditions than the parameters, and hence (\ref
{e33}) is a problem of many moment conditions with increasing number
of moments studied by \citet{HanPhi06}.

It is straightforward to verify that
\[
V_0\,{\equiv}\,\operatorname{Var}(m_n(g_0, X))\,{=}\,\operatorname{diag}\bigl\{E\bigl(\rho
(Z,g_0)^2I_{(W\in
R_{1}^n)}\bigr),\ldots,E\bigl(\rho(Z,g_0)^2I_{(W\in
R_{k_n^d}^n)}\bigr)\bigr\}.\vadjust{\goodbreak}
\]
For each $g\in\mathcal{H}$, and $j=1,\ldots,k_n^d$, write $\bar
{m}_{nj}(g)=\frac{1}{n}\sum_{i=1}^nm_{nj}(g, X_i)$
and $\bar{m}_n(g)=(\bar{m}_{n1}(g),\ldots,\bar{m}_{nk_n^d}(g))^T$. Instead
of $g_0$, we\vspace*{2pt} construct the posterior for its approximating function
inside $\mathcal{H}_n$. Under some regularity conditions, for each
fixed $k$, $\bar{m}_n(g_0)$ would satisfy the central limit theorem:
for any $\alpha\in\mathbb{R}^{k}$, as $n$ goes to infinity,
%
%
\begin{equation}\label{e34}
\Biggl|P\bigl(\sqrt{n}V_0^{-1/2}\bar{m}_n(g_0)\leq\alpha\bigr)-\prod
_{i=1}^{k}\Phi
(\alpha_i)\Biggr|\rightarrow0.
\end{equation}

This motivates a likelihood function on the sieve space $\mathcal{H}_n$,
\[
\operatorname{LIL}(g_b)\propto\exp\biggl(-\frac{n}{2}\bar
{m}_n(g_b)^TV_0^{-1}\bar{m}_n(g_b)\biggr).
\]
According to \citet{Kim02}, the function LIL$(g_b)$ can be more
appropriately interpreted as the best approximation to the true
likelihood function under the conditional moment restriction by
minimizing the Kullback--Leibler divergence, which is known as the
\textit{limited information likelihood} (LIL). Note that LIL$(g_b)$ is
not feasible, as $V_0$ depends on the unknown function $g_0$; therefore
\citet{Kim02} suggested replacing $V_0$ with a constant matrix (not
dependent on $g_0$), while maintaining the order of each element. For
each element on the diagonal, suppose we have the integration mean
value theorem: for some $w^*\in R_j^n$,
\[
E\bigl(\rho(Z,g_0)^2I_{(W\in R_{j}^n)}\bigr)=E\bigl(\rho(Z,g_0)^2|W=w^*\bigr)P(W\in
R_j^n)=O\bigl(P(W\in R_j^n)\bigr)
\]
provided that $\sup_{w\in[0,1]^d} E[\rho(Z,g_0)^2|w]<\infty$. Hence
each diagonal element of $V_0$ is of the same order as $P(W\in R_j^n)$.
We replace $V_0$ by
\[
\hat{V}=\operatorname{diag}\{\hat{v}_1,\ldots,\hat{v}_{k_n^d}\}\qquad
\mbox{where }
\hat{v}_j=\frac{1}{n}\sum_{i=1}^nI_{(W_i\in R_j^n)}.
\]
Each $\hat{v}_j$ is a consistent estimate of $P(W\in R_j^n)$. We thus
obtain the feasible LIL to be used as the likelihood function
throughout this paper,
%
%
\begin{equation}\label{e35}
L(g_b)=\exp\biggl(-\frac{n}{2}\bar{m}_n(g_b)^T\hat{V}^{-1}\bar
{m}_n(g_b)\biggr).
\end{equation}

The feasible likelihood puts more weights on the moment conditions with
smaller variance, having the same spirit of the optimal weight matrix
in \textit{generalized method of moments} [\citet{Han82}]. A more
refined approach can be based on a second-stage estimation
of $V_0$, where a consistent first-stage estimator of $g_0$ is used if
$g_0$ is assumed to be point identified.
However, it turns out that $V_0$ does not have to be estimated very
precisely in order to achieve the posterior consistency for the inference
on $g$. We will show that our simple estimator ${\hat V}$ is already
good enough for proving posterior consistency in the development to
be described below and is simple for practical computations.\vadjust{\goodbreak}

For the approximated Gaussian likelihood function (\ref{e35}), the
sample risk functional defined in Section \ref{sec2} is given by
%
%
\begin{equation}\label{gbar}
\bar{G}(g_b)\equiv\bar{m}_n(g_b)^T\hat{V}^{-1}\bar{m}_n(g_b).
\end{equation}
Let
\[
\mathcal{F}_{n}=\Biggl\{\sum_{i=1}^{q_n}b_i\phi_i(x)\dvtx\max_{i\leq
q_n}|b_i|\leq B_n\Biggr\}
\]
for some\vspace*{1pt} sequence $B_n\rightarrow\infty$; then we partition the sieve
space into $\mathcal{H}_n=\mathcal{F}_n\cup\mathcal{F}_n^c$. Under some
regularity conditions, it can be shown that\footnote{We will verify
this for the nonparametric IV regression model in Section \ref{sec4}.}
$\bar{G}$
converges in probability to the risk functional
%
%
\begin{equation}\label{G}\qquad
G(g)=E_W\{[E(\rho(Z,g)|W)]^2\}=\int_{[0,1]^d}
\bigl[E\bigl(\rho(Z,g)|W=w\bigr)\bigr]^2\,dF_W(w)
\end{equation}
uniformly on $\mathcal{F}_n$.

\subsection{Identification and ill-posedness}\label{sec32}
The identification of $g_0$ is characterized by minimizing $G$. To be
specific, define the identified region for $g_0$,
\[
\Theta_I=\bigl\{g\in\mathcal{H}\dvtx E\bigl(\rho(Z, g)|W=w\bigr)=0\mbox{ for almost
all }
w\in[0,1]^d\bigr\},
\]
which is assumed to be nonempty, then
\[
\Theta_I=\mathop{\arg\min}_{g\in\mathcal{H}}G(g)=\{g\in\mathcal{H}\dvtx
G(g)=0\}.
\]
If $\Theta_I$ is a singleton, then $\Theta_I=\{g_0\}$. Otherwise $g_0$
is partially identified on~$\Theta_I$; see, for example, \citet{San}.

In the conditional moment restriction literature, the problem of
identification and estimation of $g_0$ is well known to be \textit{ill
posed}. The ill-posed problem was postulated in detail by Kress
[(\citeyear{Kre99}),
Chapter 15], which occurs, in our context, if one of the following
three properties does not hold: (1) there exist solutions to $G(g)=0$,
and here we assume $g_0\in\Theta_I$; (2) the solution is unique, that
is, $\Theta_I$ is a singleton; (3) the solution is continuously
dependent on the data; that is, roughly speaking, when $G(g)$ is close
to zero, $g$ should be close to $\Theta_I$. However, when $g_0$ depends
on the endogenous variable $X$, the third property may fail because for
any $\varepsilon>0$, there are sequences $\{g_n\}_{n=1}^{\infty}\subset
\mathcal{H}$ such that
\[
\liminf_{n\rightarrow\infty}\inf_{g_n\notin\Theta_I^{\varepsilon}}G(g_n)=0.
\]
Throughout this paper, we call such a problem as the \textit{type-III
ill-posed inverse problem}. In order to achieve the posterior
consistency, we need certain regularization scheme to make the metric
$d(g,\Theta_I)$ be continuous with respect to the risk functional $G(g)$.\vadjust{\goodbreak}

While the literature puts a primary interest on dealing with the
type-III ill-posedness [\citet{HalHor05}, etc.], there are
relatively fewer results that deal with the second type of
ill-posedness: $\Theta_I$ is not necessarily a singleton. In this
paper, we also allow $g_0$ to be only partially identified\footnote{In
this paper, the \textit{partial identification} is meant in the
frequentist sense, as opposed to the Bayesian identification. See a
recent work by \citet{FloSimN3} for a discussion of these
concepts.} by the conditional moment restriction (\ref{e31}). Such a
treatment arises for two reasons. First, when the conditional moment
restriction is given by the nonparametric instrumental variable
regression (Example \ref{examp13}), the identification of $g_0$ depends
on the
completeness of the conditional distribution of $X|W$ [\citet
{NewPow03}]; however, the completeness assumption is hard to verify if the
conditional distribution of $X|W$ does not belong to the exponential
family. \citet{SevTri06} explored identification issues
with these models and noted that the point identification can easily
fail; see Example 3.2 of \citet{SevTri06}. For another
reason, sometimes instead of $g_0$ itself, we are only interested in a
particular characteristic of it, for example, its linear functional
$h(g_0)$. For example, in the nonparametric IV regression, if $g_0(x)$
represents the inverse demand function, then its consumer surplus at
some level $x^*$ can be written as a~functional $h(g_0)=\int
_0^{x^*}g_0(x)\,dx-g_0(x^*)x^*$. In this case, the identification of~$g_0$
might not be necessary; as \citet{SevTri06} showed,
even if~$g_0$ is not identified, it is still possible to point identify
its functional~$h(g_0)$.\vspace*{-2pt}


\subsection{Prior specification} We will apply Theorems \ref{t21} and
\ref{t22} to three types of
priors: (i) truncated prior, (ii) thin-tail prior and (iii) normal
prior. In this section we will focus on the first two types of priors,
with which more generally consistent results can be
derived.\footnote{We will describe the normal prior in a later section (Section
\ref{sec44})
since the technique
used is somewhat different, which handles mainly the situation of the
NPIV model in an identifiable situation.}

\vspace*{6pt}
\textit{Truncated prior}. The prior is supported only on $\mathcal{F}_n$.
In particular, we consider the uniform and truncated normal priors,
respectively,
\begin{eqnarray*} \mbox{uniform prior } \pi(b)&=&\prod
_{i=1}^{q_n}I(|b_i|\leq B_n);\\[-2pt]
\mbox{truncated normal } \pi(b)&=&\prod_{i=1}^{q_n}\frac
{f(b_i)I(|b_i|\leq B_n)}{P(|Z_i|\leq B_n)},
\end{eqnarray*}
where $\{Z_i\}_{i=1}^{q_n}$ are i.i.d. random variables from
$N(0,\sigma
^2)$ for some $\sigma^2>0$, and $f(\cdot)$ is the probability density
function of $Z_i$. The tail probability
\[
\pi(g_b\in\mathcal{F}_n^c)=0.\vspace*{-2pt}
\]

\textit{Thin-tail prior}.
The prior $\pi$ on $b\in\mathbb{R}^{q_n}$ is defined such that the
density is
symmetric in all directions, and $\|b\|^r$ follows\vadjust{\goodbreak} an exponential
distribution with mean $\beta^{-r}$ (for some $\beta>0$, $r>0$). Here
$\|b\|$ denotes a Euclidean norm,
\[
\pi(\|b\|^r>u^r)=e^{-\beta^ru^r},
\]
which, together with the spherical symmetry, is enough to derive the
density function,
%
%
\begin{equation}\label{e38}
\pi(b) =\frac{r\| b\|^{r-q_n}\beta^r e^{-\beta^r\| b\|^r }}{S_{q_n}},
\end{equation}
where $S_{q_n}$ is the area of the $q_n-1$-dimensional unit sphere in
Euclidean norm. For this prior, the parameter $1/\beta$ is roughly the
radius of most
of the prior mass, and $r$ denotes the thinness of the tails outside.
The bigger the $r$ is, the thinner the tail.

This prior is very similar to the class of distributions defined in
\citet{Azz86}. Both allow any positive power of the distance to the
origin to be
placed on the exponent. Our density is slightly different and does not,
in general, include the normal density exactly. However, it is derived
in a way so that the tail probability has an exact expression. Hence it
is convenient to impose a regularity condition on the tail probability.
\vspace*{8pt}

Florens and Simoni (\citeyear{FloSimN1}, \citeyear{FloSimN2})
placed a Gaussian prior whose
variance decreases to zero with the sample size. Our priors specified
here are similar to theirs in the sense that the prior tail probability
is small: when the truncated prior is used, $\pi(g_b\in\mathcal
{F}_n^c)=0$; when the thin-tail prior is used, $\pi(g_b\in\mathcal
{F}_n^c)$ decreases exponentially fast in $n$. Both types of priors
ensure that
\[
P\bigl(G(g_b)\geq\delta_n|X^n\bigr)=o_p(1)
\]
for some decaying sequence $\delta_n>0$ that depends on the convergence
rate of ${\sup_{\mathcal{F}_n}}|\bar{G}(g)-G(g)|$. The technique of using
a prior that decays exponentially fast outside a bounded sieve set is
commonly used in the nonparametric posterior consistency literature;
see, for example, \citet{GhoRam03}, \citet{GhoRoy06},
\citet{ChoSch07}, \citet{Wal03} and many references therein.

However, there is an important difference between Florens and Simoni's
prior settings (2009a) and our own. While \citet{FloSimN1}
put their prior on an infinite-dimensional function space, they require
the variance of the Gaussian prior to shrink to zero as a
regularization scheme in order to achieve the posterior consistency. In
contrast, our prior is placed directly on the sieve coefficients
$(b_1,\ldots,b_{q_n})$ in a finite-dimensional vector space, and neither
the truncated prior nor the thin-tail prior shrinks to a~point mass.\vadjust{\goodbreak}
When $q_n$ grows slowly with $n$, it can be shown that\footnote{We
will verify this for the nonparametric IV regression model.} for any
$\varepsilon>0$,
\[
\inf_{g_b\in\mathcal{H}_n, d(g_b,\Theta_I)\geq\varepsilon
}G(g_b)\succ
\delta_n;
\]
hence the distinguishing ability condition in Theorem \ref{t22} is
satisfied. As a~result, in our procedure it is the fact that $q_n$
grows slowly that plays the role of regularization instead of a
shrinking prior. Later in Section \ref{sec44}, we will also verify that
with a
suitably chosen $q_n$, a nonshrinking normal prior can be used to
achieve the posterior consistency in the identified NPIV model.

\subsection{Posterior consistency}

The following assumptions are imposed.
%
%
\begin{assum}\label{iid}
The data $X^n=(X_1,\ldots,X_n)$ are independent and identically distributed.
\end{assum}
%
%
\begin{assum} \label{a32} There exists a positive sequence $\lambda
_n\rightarrow0$ such that
\[
{\sup_{g\in\mathcal{F}_n}}|\bar{G}(g)-G(g)|=O_p(\lambda_n).
\]
\end{assum}

Since $\mathcal{F}_n$ is compact in $\mathcal{H}_n$, as long as the
radius of $\mathcal{F}_n$ grows slowly, the uniform convergence
condition in Assumption \ref{a32} can be shown using similar
techniques to those in \citet{HanPhi06}. We will verify it
for the nonparametric IV regression example in Section \ref{sec4}.
%
%
\begin{assum}\label{a33} (i) $\{\phi_1,\phi_2,\ldots,\phi_{q_n}\}$ forms
an orthonormal basis of $\mathcal{H}_n$ such that $E(\phi_i(X)\phi
_j(X))=\delta_{ij}$, the Kronecker $\delta$.

(ii) There exist $g_0\in\Theta_I$, and $g_{q_n}^*=\sum
_{i=1}^{q_n}b_i^*\phi_i\in\mathcal{H}_n$ such that $\|g_{q_n}^*-g_0\|
_s=o(1)$ as $q_n\rightarrow\infty$.

\end{assum}

The existence of $g_{q_n}^*$ is simply implied by the definition of a
sieve space. It is satisfied by the spaces that are spanned by commonly
used sieve basis functions such as splines, power series, wavelets and
Fourier series. For example, if the parameter\vspace*{-1pt} space is a Sobolev space
$\mathcal{W}^2_{p}[0,1]^{d_x}$, where $d_x=\dim(X)$, and \mbox{$\|\cdot\|_s$} is
the Sobolev norm,\vspace*{1pt} then $\|g_{q_n}^*-g_0\|_s=O(q_n^{-p/d_x})$ for some
$p>0$; see, for example, Kress [(\citeyear{Kre99}), Chapter 8] and
\citet{Che}; see
also \citet{Sch81} and \citet{Mey90} for splines and orthogonal
wavelets in other function spaces.
%
%
\begin{assum} \label{a34} There exists $C>0$ such that $\forall g_1,
g_2\in\mathcal{H}$,
\[
E|\rho(Z,g_1)-\rho(Z, g_2)|\leq CE|g_1(X)-g_2(X)|.
\]
\end{assum}

This assumption is trivially satisfied by the nonparametric IV
regression in Example \ref{examp13}. Here we give another example that satisfies
this assumption.

%
\begin{exm}[(Nonparametric quantile IV regression)] \label{ex32}
Consider the~mo\-del in Example \ref{examp14}, in which the conditional moment
restriction is given by
\[
E(\rho(Z,g_0)|W, g_0)=0,\qquad \rho(Z,g_0)=I_{(y\leq
g_0(X))}-\gamma.
\]
It is straightforward to verify that for any $g_1, g_2$,
\begin{eqnarray*}
E|\rho(Z,g_1)-\rho(Z, g_2)|&=& E\bigl|I_{(g_1(X)\leq y\leq
g_2(X))}+I_{(g_2(X)\leq y\leq g_1(X))}\bigr|\\
&=& E\bigl[P\bigl(g_1(X)\leq y\leq g_2(X)|X\bigr)\bigr]\\
&&{}+E\bigl[P\bigl(g_2(X)\leq y\leq g_1(X)|X\bigr)\bigr].
\end{eqnarray*}
Suppose there exists a constant $C>0$ such that $F_{y|X}(\cdot)$, the
conditional c.d.f. of $y|X$, satisfies
\[
|F_{y|x}(y_1)-F_{y|x}(y_2)|\leq C|y_1-y_2|
\]
for any $y_1, y_2\in\mathbb{R}$ and $x$ in the support of $X$. Then the
first term on the right-hand side is bounded by
\begin{eqnarray*}
E\bigl[P\bigl(g_1(X)\leq y\leq g_2(X)|X\bigr)\bigr]&\leq&
E|F_{y|X}(g_2(X))-F_{y|X}(g_1(X))|\\
&\leq& CE|g_2(X)-g_1(X)|.
\end{eqnarray*}
Likewise, $E[P(g_2(X)\leq y\leq g_1(X)|X)]\leq CE|g_2(X)-g_1(X)|$.
Therefore Assumption \ref{a34} is satisfied.
\end{exm}

Define
%
%
\begin{equation} \label{e39}
\gamma_n={\sup_{g\in\mathcal{F}_n, w\in[0,1]^d}}
\bigl|E\bigl(\rho(Z, g)|W=w\bigr)\bigr|+1.
\end{equation}

We are able to verify the conditions in Theorem \ref{t21} with the
previous assumptions, and establish the following theorem:
%
%
\begin{theorem}[(Risk consistency: truncated prior)] \label{t31} Suppose
$q_n=o(n)$ and $B_n=o(n)$. Assume $\delta_n=O(1)$ is such that there
exists $g_0\in\Theta_I$ whose sieve approximation $g_{q_n}^*$ satisfies
\[
\max\biggl\{G(g_{q_n}^*), \lambda_n, \frac{q_n}{n}\log(\gamma_nn)\biggr\}
=o(\delta_n).
\]
Then when either the uniform prior or the truncated normal prior is
used, under Assumptions \ref{iid}--\ref{a34},
\[
P\bigl(G(g_b)<\delta_n|X^n\bigr)\rightarrow^p1.
\]
\end{theorem}

In the following theorem, write $\lambda(B_n)=\lambda_n$ and $\gamma
(B_n)=\gamma_n$ to indicate the dependence of $\lambda_n$ and $\gamma
_n$ on $B_n$, defined in Assumption \ref{a32} and (\ref{e39}),
respectively.
%
%
\begin{theorem}[(Risk consistency: thin-tail prior)] \label{t32} Suppose
there exists $g_0\in\Theta_I$ with $g_{q_n}^*$ being its sieve
approximation in $\mathcal{H}_n$, and a sequence $B_n^*\rightarrow
\infty$ such that $\max\{G(g_{q_n}^*), \lambda(B_n^*), \gamma
(B_n^*)e^{-n\lambda(B_n^*)/q_n}\}=o( B_n^{*r}/n)$. In
addition, suppose $\delta_n=O(1)$ is such that
\[
\max\bigl\{G(g_{q_n}^*), \lambda(B_n^*), \gamma
(B_n^*)e^{-n\lambda(B_n^*)/q_n}\bigr\}=o(\delta_n).
\]
Then under Assumptions \ref{iid}--\ref{a34},
\[
P\bigl(G(g_b)<\delta_n|X^n\bigr)\rightarrow^p1.
\]
\end{theorem}
%
%
\begin{remark}
(1) We will show in the next section that in the nonparametric
IV regression model, $\gamma_n=O(q_nB_n)$. For the nonparametric
quantile IV regression in Example \ref{ex32}, $\gamma_n$ is a constant
that is bounded away from zero.

(2) Under the conditions of Theorems \ref{t31} and \ref{t32},
$\delta_n$ can be fixed as a~constant. Namely, $\forall\delta>0$,
\[
P\bigl(G(g_b)>\delta|X^n\bigr)=o_p(1).
\]
Roughly speaking, the posterior distribution is asymptotically
supported on the set where $G$ is minimized. This result has many
important applications. For example, in the binary treatment effect
study, let $Y\in\{0,1\}$ indicate whether a treatment is successful,
which is associated with a covariate $X$. Suppose we model the success
probability $P(Y=1|X=x)$ by a nonparametric function $g(x)$. In this model,
\[
G(g)=E_X\{[E(Y|X)-g(X)]^2\}=\bigl\|P(Y=1|X)-g(X)\bigr\|_s^2,
\]
where $\|g\|_s^2=E(g(X)^2)$. By Theorems \ref{t31}, \ref{t32}, for
any $\varepsilon>0$, the posterior
\[
P\bigl(\bigl\|P(Y=1|X)-g_b(X)\bigr\|_s^2<\varepsilon|\mbox{Data}\bigr)\rightarrow^p1,
\]
which implies that the posterior of $g_b$ can recover the success
probability arbitrarily well with high probability.

(3) In data mining, this type of result is sometimes
called the ``risk consistency.'' For example, if $G$ was the
classification risk, the risk consistency result would show that the
posterior would effectively minimize the misclassification error.
The current definition of $G$, however, is not the classification
risk. In nonparametric regression and in the NPIV example, the risk
$G$ becomes, respectively, $E_W\{[E(Y|W)-g(W)]^2\}$ and
$E_W\{[E(Y|W)-E(g(X)|W)]^2\}$, which is related to how much $E(Y|W)$
would be missed if it was estimated by (something derived from) $g$.
\end{remark}

The following two theorems establish the posterior consistency without
assuming the compactness of the parameter space $\mathcal{H}$.\vadjust{\goodbreak}
%
%
\begin{theorem}[(Posterior consistency: truncated prior)] \label{t33}
Suppose there exists $g_0\in\Theta_I$ whose sieve approximation
$g_{q_n}^*$ satisfies $\forall\varepsilon>0$
%
%
\begin{equation} \label{e310}
\max\biggl\{G(g_{q_n}^*), \lambda_n, \frac{q_n}{n}\log(\gamma_nn)\biggr\}
=o\Bigl(\inf
_{g\in\mathcal{H}_n, g\notin\Theta_I^{\varepsilon}}G(g)\Bigr).
\end{equation}
Then under Assumptions \ref{iid}--\ref{a34}, for any $\varepsilon>0$,
\[
P\bigl(d(g_b,\Theta_I)<\varepsilon|X^n\bigr)\rightarrow^p1.
\]
\end{theorem}
%
%
\begin{theorem}[(Posterior consistency: thin-tail prior)] \label{t34}
Suppose there exists $g_0\in\Theta_I$ with $g_{q_n}^*$ being its sieve
approximation in $\mathcal{H}_n$, and a sequence $B_n^*\rightarrow
\infty$ such that $\max\{G(g_{q_n}^*), \lambda(B_n^*), \gamma
(B_n^*)e^{-n\lambda(B_n^*)/q_n}\}=o( B_n^{*r}/n)$. In
addition, suppose $\forall\varepsilon>0$,
%
%
\begin{equation} \label{e311}
\max\bigl\{G(g_{q_n}^*), \lambda(B_n^*),
\gamma
(B_n^*)e^{-n\lambda(B_n^*)/q_n}\bigr\}=o\Bigl(\inf_{g\in\mathcal{H}_n,
g\notin
\Theta_I^{\varepsilon}}G(g)\Bigr).
\end{equation}
Then under Assumptions \ref{iid}--\ref{a34}, for any $\varepsilon>0$,
\[
P\bigl(d(g_b,\Theta_I)<\varepsilon|X^n\bigr)\rightarrow^p1.
\]
\end{theorem}
%
%
\begin{remark}
$\!\!$(1) The restriction $\lambda(B_n^*)=o(B_n^{*r}/n)$ in both
Theorems~\ref{t32} and~\ref{t34} requires that $r$, the thin-tail
prior parameter, should not be too small; otherwise, no such $B_n^*$
exists. In the NPIV model which will be illustrated in the next
section, we need $r>6d+4$, where $d=\dim(W)$.

(2) Conditions (\ref{e310}) and (\ref{e311}) are similar to
Chen and Pouzo's [(\citeyear{ChePou}),\allowbreak condition (3.1)], where they
require that $q_n$
grow slowly enough so that
$\inf_{g\in\mathcal{H}_n, g\notin\Theta
_I^{\varepsilon}}G(g)$ does not decrease too fast for any fixed
$\varepsilon
>0$. This will also be illustrated in Section \ref{sec4}.
\end{remark}

Let $h(g_0)$ be a linear functional of $g_0$, whose practical meaning
may be of direct interest. For example, if $h(g_0)=E[g_0(X)\omega(X)]$
for some weight function $\omega$, then with proper choices of $\omega
$, $h$ can be used to test some special properties of $g_0$, such as
the monotonicity, the convexity, etc. \citet{San11}. On the other
hand, $h$ itself may have interesting meanings. For example, when $g_0$
denotes the inverse demand function in nonparametric regression,~$h(g_0)$ can be the consumer surplus [\citet{San}].
\citet{SevTri06} have provided conditions to point identify $h(g_0)$
even if $g_0$ itself is not identified.
%
%
\begin{exm}\label{examp32}
Suppose we want to test whether the unknown function $g_0$ is weakly
increasing. Note that any weakly increasing function $g(x)$ must
satisfy $\int_{-\pi}^{\pi}\sin(x)g(x)\,dx\geq0$; hence the functional
of interest here is $h(g_0)=\int_{-\pi}^{\pi}\sin(x)g_0(x)\,dx$. Suppose
the joint distribution of $(X,W)$ has density function $f_{XW}(x,w)$.
By \citet{SevTri06}, $h(g_0)$ is point identified, if there
exists $p(w)$ such that $E[p(W)^2]<\infty$ and $E(p(W)|X)=\sin
(X)/f_X(X)$ almost surely.\vadjust{\goodbreak}
\end{exm}

Theorems \ref{t33} and \ref{t34} imply a flexible way to
consistently estimate $h$ without identifying $g_0$. In the following
assumption, condition (i) assumes the point identification of $h(g_0)$.
Condition (ii) requires the uniform continuity of $h$, which is
satisfied when $h(g)=E[g(X)\omega(X)]$ if $\sup_{x}|w(x)|<\infty$ and
$E|g_1-g_2|\leq C\|g_1(X)-g_2(X)\|_s$ for any $g_1, g_2\in\mathcal{H}$.
%
%
\begin{assum}\label{h}
(i) $\{h(g)\dvtx g\in\Theta_I\}=\{h(g_0)\}$; (ii) $h\dvtx(\mathcal{H},
\mbox{$\|\cdot\|_s$})\rightarrow\mathbb{R}$ is uniformly continuous.
\end{assum}
%
%
\begin{cor} \label{cor1} Suppose the assumptions of Theorem \ref{t33}
(if the truncated priors are used) and Theorem \ref{t34} (if the
thin-tail prior is used) are satisfied. In addition, suppose Assumption
\ref{h} holds. When $g_0$ is not necessarily point identified,
$\forall
\delta>0$,
\[
P\bigl(|h(g_b)-h(g_0)|<\delta|X^n\bigr)\rightarrow^p1.
\]
\end{cor}

\section{Nonparametric instrumental variable regression}\label{sec4}
\subsection{The model}
$\!\!\!$The nonparametric instrumental variable regression (NPIV) model is
given by
\[
Y=g_0(X)+\varepsilon,
\]
where $X$ is endogenous, which is correlated with $\varepsilon$. We
consider the following parameter space and the norm \mbox{$\|\cdot\|_s$}:
\[
\mathcal{H}=L^2(X)=\{g\dvtx E[g(X)^2]<\infty\},\qquad
\|g\|_s^2=E[g(X)^2].
\]
In addition, suppose we observe an instrumental variable $W\in[0,1]^d$
such that $E(\varepsilon|W)=0$. Applications of instrumental variables can
be found in many standard econometrics texts, for example, \citet
{Han}. Let $Z=(Y,X)$; the NPIV model is then essentially a conditional
moment restricted model with $\rho(Z,g)=Y-g(X)$.

Let $\{\phi_1,\phi_2,\ldots\}$ be a set of orthonormal basis functions of
$L^2(X)$. We consider the sieve space $\mathcal{H}_n=\{g\in L^2(X)\dvtx
g=\sum_{i=1}^{q_n}b_i\phi_i\}$, which can be partitioned into
$\mathcal
{H}_n=\mathcal{F}_n\cup\mathcal{F}_n^c$, where $\mathcal{F}_n=\{
\sum
_{i=1}^{q_n}b_i\phi_i\in\mathcal{H}_n, {\max_{i\leq q_n}}|b_i|\leq
B_n \}
$ as in\vspace*{1pt} Section \ref{sec3}.

We apply the feasible LIL (\ref{e35}) to construct the posterior. The
log-likelihood involves the sample risk functional
\[
\bar{G}(g)=\sum_{j=1}^{k_n^d}\Biggl(\frac{1}{n}\sum
_{i=1}^n\bigl(Y_i-g(X_i)\bigr)I_{(W_i\in R_j^n)}\Biggr)^2\hat{v}_j^{-1},
\]
which later will be shown to uniformly converge to
\[
G(g)=E_W\bigl\{\bigl[E\bigl(Y-g(X)|W\bigr)\bigr]^2\bigr\}
\]
over $\mathcal{F}_n$. The identified region $\Theta_I$ is defined as a
subset of $L^2(X)$ on which \mbox{$G(g)=0$}.

\subsection{Risk consistency}\label{sec42}
Under mild conditions, we can derive the convergence rate of ${\sup
_{g\in
\mathcal{F}_n}|}\bar{G}(g)-G(g)|$. The following assumptions are imposed.

%
\begin{assum} \label{a41} (i) $k_n^{-d}=O(\min_{j\leq k_n^d}P(W\in
R_j^n))$;

(ii) $\max_{j\leq k_n^d}P(W\in R_j^n)=O(k_n^{-d})$.
\end{assum}

This assumption is satisfied, for example, when $W$ has a continuous
density function on $[0,1]^d$ that is bounded away from both zero and infinity.
%
%
\begin{assum} \label{a42} There exists $C>0$ such that for all
$i=1,\ldots,q_n$:

\begin{longlist}
\item
$\sup_wE(Y^2|W=w)<C$, $\sup_wE(\phi_i(X)^2|W=w)<C$;

\item
$E(Y|W=w)$ is Lipschitz continuous with respect to $w$ on
$[0,1]^d$;

\item for any $w_1, w_2\in[0,1]^d$,
\[
\bigl|E\bigl(\phi_i(X)|W=w_1\bigr)-E\bigl(\phi_i(X)|W=w_2\bigr)\bigr|
\leq C\|w_1-w_2\|.
\]
\end{longlist}
\end{assum}

Condition (iii) requires that the family $\{E(\phi_i(X)|W=w)\dvtx i\leq
q_n\}$ is Lipschitz equicontinuous on $[0,1]^d$, which is satisfied,
for example, when $X$ has a density function that is bounded away from
zero on the support of $X$; in addition, $X|W$ has a conditional
density function $f_{X|W}$ such that for some $C>0$,
\[
\bigl|f_{X|W}(x|w_1)-f_{X|W}(x|w_2)\bigr|\leq C\|w_1-w_2\|
\]
for all\vadjust{\goodbreak} $x$ and $w_1, w_2\in[0,1]^d$.\footnote{This is simple to show:
for any $w_1, w_2$,
$|E(\phi_i(X)|W=w_1)-E(\phi_i(X)|W=w_2)|\leq(\inf f_X(x))^{-1}\int
|\phi_i(x)f_X(x)||f_{X|W}(x|w_1)-f_{X|W}(x|w_2)|\,dx\leq C\|w_1-w_2\|
E|\phi_i(X)|\leq C'\|w_1-w_2\|$, where the fact that $E|\phi_i(X)|$ is
bounded away from infinity is guaranteed by condition (i).}
%
%
\begin{assum} \label{a43} There exist $g_0\in\Theta_I$,
$g_{q_n}^*=\sum
_{i=1}^{q_n}b_i^*\phi_i$ with $\sum_{i=1}^{\infty}b_i^{*2}<\infty$, and
a positive sequence $\{\eta_{j}\}_{j=1}^{\infty}$ that strictly
decreases to zero as $j\rightarrow\infty$ such that $\|g_{q_n}^*-g_0\|
_s=O(\eta_{q_n})$ as $q_n\rightarrow\infty$. (We will choose
$g^*_{q_n}$ to be the projection of
$g_0$ onto ${\mathcal H}_n$, unless otherwise noted.)
\end{assum}

Examples of the rate $\eta_{q_n}$ are discussed earlier behind
Assumption \ref{a33}.
%
%
\begin{theorem}\label{t41} Assume $q_n^2B_n^2=o(\min\{\sqrt{n}/k_n^{3d/2},
k_n\})$. Then under Assumptions \ref{iid}, \ref{a41}, \ref{a42},
\[
{\sup_{g\in\mathcal{F}_n}}|\bar{G}(g)-G(g)|=O_p\biggl(\frac
{q_n^2B_n^2k_n^{3d/2}}{\sqrt{n}}+\frac{q_n^2B_n^2}{k_n}\biggr).
\]
\end{theorem}

Define a semi-norm \mbox{$\|\cdot\|_w$}, which is weaker than \mbox{$\|\cdot\|_s$}, as
%
%
\begin{equation} 
\|g\|_w^2=E_W\{(E(g(X)|W))^2\}.
\end{equation}
It can be easily verified that $\|\cdot\|_w$ satisfies the triangular
inequality, but $\|g\|_w=0$ does not necessarily imply $g=0$ if the
conditional distribution $X|W$ is not complete. Note that $G(g)=\|
g_0-g\|_w^2$; hence this semi-norm induces an equivalence class
characterized by the identified region $\Theta_I=\{g\in L^2(X)\dvtx
E(Y-g(X)|W)=0\mbox{, a.s.}\} $, such that $\|g-g_0\|_w=0$ if and only if $g\in
\Theta_I$. In other words, we can say that $g_0$ \textit{is weakly
identified} under $\|\cdot\|_w$, since for any $g\in\Theta_I$, $g$ and
$g_0$ are equivalent under $\|\cdot\|_w$.

The following theorem is a straightforward application of Theorems \ref
{t31} and~\ref{t32}:
%
%
\begin{theorem}[(Risk-consistency)] \label{t42} Under Assumptions \ref{iid},
\ref{a41}--\ref{a43}, suppose $\delta_n=O(1)$ is such that:

\begin{longlist}
\item
for the truncated priors assuming $q_n^2B_n^2=o(n^{1/(3d+2)})$,
\[
\max\biggl\{\eta_{q_n}^2, q_n^2B_n^2\biggl(\frac{k_n^{3d/2}}{\sqrt
{n}}+\frac{1}{k_n}\biggr)\biggr\}=o(\delta_n),
\]
\item for the thin-tail prior with $r>6d+4$, assuming $q_n=o(n^{1/(6d+4)-1/r})$,
\[
\max\biggl\{\eta_{q_n}^2, n^{2/(r-2)}q_n^{2r/(r-2)}\biggl(\frac
{k_n^{3d/2}}{\sqrt{n}}+\frac{1}{k_n}\biggr)^{r/(r-2)}\biggr\}
=o(\delta_n),
\]
then
\[
P(\|g_b-g_0\|_w>\delta_n|X^n)=o_p(1).
\]
\end{longlist}
\end{theorem}

\subsection{Ill-posedness and posterior consistency}\label{sec43}

Define
\[
T\dvtx L^2(X)\rightarrow\{\zeta\dvtx E[\zeta(W)^2]<\infty\},\qquad
T(g)=E(g(X)|W)
\]
and
write $E(Y|W=w)\equiv\zeta(w)$. Then the NPIV model can be
equivalently written as
%
%
\begin{equation}\label{e41}
Tg_0=\zeta.
\end{equation}

Under Assumption \ref{a44}, $T$ is a compact linear operator [see
\citet{CarFloRen}], and therefore is continuous. Equation (\ref
{e41}) is usually called the \textit{Fredholm integral equation of the
first kind}.
%
%
\begin{assum} \label{a44} The joint distribution $(Y,X,W)$ is
absolutely continuous with respect to the Lebesgue measure. In
addition, suppose $f_{XW}(x,w)$, $f_X(x)$, $f_W(w)$ denote the density
functions of $(X,W)$, $X$ and $W$, respectively, then
\[
\iint\biggl(\frac{f_{XW}(x,w)}{f_X(x)f_W(w)}
\biggr)^2f_X(x)f_W(w)\,dx\,dw<\infty.
\]
\end{assum}

As described before, the problem of inference about $g_0$ is ill-posed
in two aspects. The first ill-posedness comes from the identification,
which depends on the invertibility of $T$. If $T$ is nonsingular, in
which case its null space is~$\{0\}$, $g_0$ can be point identified by
$g_0=T^{-1}\zeta$, but not otherwise. See \citet{SevTri06}
and \citet{DHa11} for detailed descriptions of the
identification issues.

Even when $g_0$ is identified, in which case $T^{-1}$ exists, as
pointed out by \citet{Fl03} and \citet{HalHor05}, since
$L^2(X)$ is of infinite dimension, and $T$ is compact, $T^{-1}$ is not
bounded (therefore is not continuous). As a result, small inaccuracy in
the estimation of $\zeta$ can lead to large inaccuracy in the
estimation of $g_0$, which is known as the type-III ill-posed inverse
problem described in Section~\ref{sec32}. When $g_0$ is partially identified,
this problem is still present when
\[
\liminf_{n\rightarrow\infty}\inf_{g\in\mathcal{H}_n,g\notin
\Theta
_I^{\varepsilon}}G(g)
=\liminf_{n\rightarrow\infty}\inf_{g\in\mathcal{H}_n,g\notin
\Theta
_I^{\varepsilon}}E\{[T(g-g_0)]^2\}=0.
\]

By Theorems \ref{t33}, \ref{t34} and \ref{t42}, in order to achieve
the posterior consistency, it suffices to verify
%
%
\begin{equation}\label{e43} \delta_n^*=o\Bigl(\inf_{g\in\mathcal{H}_n,
g\notin\Theta_I^{\varepsilon}}G(g)\Bigr),
\end{equation}
where
\begin{eqnarray*}
\mbox{for truncated prior }
\delta_n^*&=&\max\biggl\{\eta_{q_n}^2, q_n^2B_n^2\biggl(\frac
{k_n^{3d/2}}{\sqrt{n}}+\frac{1}{k_n}\biggr)\biggr\},
\\
\mbox{for thin-tail prior }
\delta_n^*&=&\max\biggl\{\eta_{q_n}^2, n^{2/(r-2)}q_n^{2r/(r-2)}
\biggl(\frac{k_n^{3d/2}}{\sqrt{n}}+\frac{1}{k_n}\biggr)^{r/(r-2)}\biggr\}.
\end{eqnarray*}
Hence it\vspace*{1pt} requires us to derive a lower bound of $\inf_{g\in\mathcal
{H}_n, g\notin\Theta_I^{\varepsilon}}G(g)$ first, and, in addition, this
lower bound should decay at a rate slower than $\delta_n^*$.

When $g_0$ is point identified and a slowly growing finite-dimensional
sieve is used, \citet{ChePou} showed the existence of such a
lower bound using the singular value decomposition of $T$. Their
approach is briefly illustrated in the following example.
%
%
\begin{exm} \label{ex41}
Let $\langle g_1,g_2\rangle_X=E[g_1(X)g_2(X)]$ denote the inner product
of two elements in $L^2(X)$, and $\{\nu_j,\phi_{1j}, \phi_{2j}\}
_{j=1}^{\infty}$ be the ordered singular value system of $T$ such that
\[
T\phi_{1j}=\nu_j\phi_{2j},\qquad \nu_1^2\geq\nu_2^2\geq
\cdots.
\]
Suppose $T$ is nonsingular, then $\{\phi_{1j}\}_{j=1}^{\infty}$ forms
an orthonormal basis of~$L^2(X)$. \citet{ChePou} showed that
when $\{\phi_{1j}\}_{j=1}^{q_n}$ is used as the basis in the
sieve\vadjust{\goodbreak}
approximation space, $\forall\varepsilon>0, \nu_{q_n}^2=O(\inf_{g\in
\mathcal{H}_n, g\notin\Theta_I^{\varepsilon}}G(g))$. Therefore, condition
(\ref{e43}) is satisfied if we assume $\delta_n^*=o(\nu_{q_n}^2)$. In
addition, suppose $\{\nu_j^2\}_{j=1}^{\infty}$ decays at a polynomial
rate $j^{-\alpha}$ for some $\alpha>0$; then we require $q_n=o(\delta
_n^{*-1/\alpha})$, a slowly growing sieve dimension.
\end{exm}

We impose the following assumption to derive a lower bound for\break $\inf
_{g\in\mathcal{H}_n, g\notin\Theta_I^{\varepsilon}}G(g)$ and verify
(\ref
{e43}), which, in the identified case, uses more general basis
functions for the sieve space. Therefore we allow the sieve basis to be
different from the eigenfunctions of $T$. A similar approach was used
by Chen and Reiss [(\citeyear{CheRei11}), Section 6.1], who used the
wavelets as the
sieve basis functions while the eigenfunctions of $T$ form a Fourier basis.

%
%
%
\begin{assum}\label{a45} There is a continuous and increasing function
$\varphi(\cdot)>0$ satisfying $\lim_{t\rightarrow0^+}\varphi(t)=0$ such
that, for $\{g_0, g_{q_n}^*, \{\eta_j\}_{j=1}^{\infty}\}$ as defined in
Assumption \ref{a43} and some constants $C_1, C_2>0$:

\begin{longlist}
\item
$\|g-g_0\|_w^2\geq C_1\sum_{j=1}^{\infty}\varphi(\eta
_j^2)|\langle
g-g_0, \phi_j\rangle_X|^2$ for all $g\in L^2(X)$;

\item $\|g_{q_n}^*-g_0\|_w^2\leq C_2\sum_j\varphi(\eta_j^2)|\langle
g_0-g_{q_n}^*, \phi_j\rangle_X|^2$.
\end{longlist}
\end{assum}
%
%
\begin{remark}
(1) This assumption implies a generalization of the relation
$\nu
_{q_n}^2=O(\inf_{g\in\mathcal{H}_n, g\notin\Theta_I^{\varepsilon}}G(g))$
in Example \ref{ex41}. In this assumption, $\{\phi_j\}_{j=1}^{\infty}$
are the basis functions whose first $q_n$ terms span the sieve
approximation space. In the identified case, $\{\phi_j\}_{j=1}^{\infty
}$ can be a general set of basis functions that is different from the
eigenfunctions of $T$. Chen and Pouzo [(\citeyear{ChePou}), Section
5.3] identified
the singular value $\nu_j^2$ of Example \ref{ex41} as a special case
of the general $\varphi(\eta_j^2)$, in which case Assumption \ref{a45}
is satisfied. In its general form, Assumption \ref{a45} is standard in
the literature for the linear ill-posed inverse problem when the
convergence rate of the estimator is studied; see, for example,
\citet{NaiPerTau05}, Chen and Pouzo [(\citeyear{ChePou}),
Assumption 5.2], Chen and Reiss
[(\citeyear{CheRei11}), Section 2.1], etc. As described above, however,
this assumption
is also needed in order to verify~(\ref{e43}) and show consistency
when general basis functions are used. \citet{BluCheKri07} provided
sufficient conditions of Assumption~\ref{a45} for the NPIV model setting.

(2) In the partially identified case when $\Theta_I$ is not a
singleton, Assumption~\ref{a45} is still satisfied, if we take $\{
\phi
_j\}_{j=1}^{\infty}$ to be the eigenfunctions of $T^*T$ that correspond
to its nonzero eigenvalues, where $T$ is the conditional expectation
operator, and $T^*$ is its adjoint. The spectral theory of compact
operators [\citet{Kre99}] implies that\vspace*{1pt} $\|T(g-g_0)\|_s^2= \sum
_{j=1}^{\infty}\nu_j^2|\langle g-g_0, \phi_j\rangle_X|^2$ for all
$g\in
L^2(X)$, where $\{\nu_j^2\}$ represent all the (nonzero) eigenvalues of~$T^*T$, and $\{\phi_j\}$ are the corresponding eigenfunctions (the zero
eigenvalues of $T^*T$ do not contribute to the right-hand side of the
spectral decomposition).
Therefore, Assumption \ref{a45} remains valid with $\varphi(\eta
_j^2)=\nu_{j}^2$,
with $\{\nu^2_{j}\}$ denoting\vspace*{1pt} the sequence of decreasing nonzero eigenvalues.
This\vadjust{\goodbreak} idea of using the spectral representation of $T^*T$ is related to
the commonly used ``general source condition'' in the literature
[\citet{Tau98} and \citet{Daretal}], where, for example,
\citet{Daretal} used this condition to derive the convergence
rate of their kernel-based Tikhonov regularized estimator in NPIV regression.

(3) When a more general sieve basis $\{\phi_j\}_{j=1}^{\infty}$
is used in the partially identified case, condition~(i) of Assumption
\ref{a45} is not generally satisfied. For example, suppose there
exists $g\in\Theta_I$, but $g\neq g_0$. By the definition of $\|\cdot\|_w$,
$\|g-g_0\|_w^2=0$, but the right-hand side of the displayed inequality
in condition~(i) is strictly positive unless $\{\phi_j\}_{j=1}^{\infty
}$ are the\vspace*{1pt} eigenfunctions of $T^*T$. To allow for more general sieve
basis in this case, a possible approach is to assume the true $g_0$ in
the data generating process to lie in a compact set~$\Theta$, for
example., a Sobolev ball [\citet{CheRei11}]. It is then not hard
to show that $\inf_{ g\in\Theta, g\notin\Theta_I^{\varepsilon}}G(g)$ is
bounded away from zero. Restricting $g_0$ inside a compact set is
actually a quite common approach in nonparametric IV regression, and
the literature is found in \citet{NewPow03},
\citet{BluCheKri07}, \citet{CheRei11}, etc. Recently,
\citet{San} extended
this approach to the partially identified case, with the compactness
restriction. We do not pursue this approach here, since our other
results on posterior consistency allow a noncompact parameter space.
%
\end{remark}

As in \citet{ChePou}, generally the degree of ill-posedness
has two types:
\begin{longlist}[(2)]
\item[(1)] \textit{mild ill-posedness}: $\varphi(\eta)=\eta
^{\alpha}$
for some $\alpha>0$.
\item[(2)] \textit{severe ill-posedness}: $\varphi(\eta)=\exp
(-\eta
^{-\alpha})$ for some $\alpha>0$.
\end{longlist}

Under Assumption \ref{a45}, it can be shown that $\varphi(\eta
_{q_n}^2)\,{=}\,O(\inf_{g\in\mathcal{H}_n, g\notin\Theta_I^{\varepsilon}}G(g))$
for any $\varepsilon>0$; see Lemma C.5 of the supplementary material.
Intuitively speaking, $\varphi(\cdot)$ is associated with the singular
values of $T$ and is related to how severe the type-III ill-posed
inverse problem is. When the nonzero singular values decay at a
polynomial rate, $\varphi$ corresponds to the mildly ill-posed case;
when the singular values decay at an exponential rate, it corresponds
to the severely ill-posed case.


Before formally presenting our posterior consistency result, we briefly
comment on the role of condition (ii) of Assumption \ref{a45}.
Assumption 5.2(ii) is the so-called ``stability condition'' in
\citet{ChePou}
that is required to hold only in terms of the sieve
approximation error on one element in $\Theta_I$. By Theorems \ref
{t33} and \ref{t34}, we require $G(g_{q_n}^*)=o(\inf_{g\in\mathcal
{H}_n, g\notin\Theta_I^{\varepsilon}}G(g))$. It can be easily shown that
$G(g_{q_n}^*)=O(\eta_{q_n}^2)$, and hence $G(g_{q_n}^*)$ was\vspace*{1pt} replaced
with $\eta_{q_n}^2$ in the condition of Theorem \ref{t42}. In addition,
condition~(i) of Assumption \ref{a45} implies that $\varphi(\eta
_{q_n}^2)=O(\inf_{g\in\mathcal{H}_n, g\notin\Theta_I^{\varepsilon}}G(g))$.
With condition~(ii) of Assumption \ref{a45}, it can be further shown that
$G(g_{q_n}^*)\,{=}\,O(\eta_{q_n}^2\varphi(\eta_{q_n}^2))$ (see Lemma C.6 in
the supplementary material). Since $\eta_{q_n}^2=o(1)$,
$G(g_{q_n}^*)=o(\varphi(\eta_{q_n}^2))=o(\inf_{g\in\mathcal{H}_n,
g\notin\Theta_I^{\varepsilon}}G(g))$ is verified.

Under this framework, we have the posterior consistency under \mbox{$\|\cdot\|_s$}:
%
%
\begin{theorem}[(Posterior consistency)] \label{t43} Under Assumptions
\ref
{iid}, \ref{a41}--\ref{a45}, suppose:

\begin{longlist}
\item
for the truncated priors assuming $q_n^2B_n^2=o(n^{1/(3d+2)})$,
%
%
\begin{equation} \label{e44}
q_n^2B_n^2\biggl(\frac{k_n^{3d/2}}{\sqrt{n}}+\frac{1}{k_n}
\biggr)=o(\varphi(\eta_{q_n}^2));
\end{equation}
\item for the thin-tail prior with $r>6d+4$, assuming $q_n=o(n^{1/(6d+4)-1/r})$,
%
%
\begin{equation} \label{e45}
n^{2/(r-2)}q_n^{2r/(r-2)}\biggl(\frac{k_n^{3d/2}}{\sqrt{n}}+\frac
{1}{k_n}\biggr)^{r/(r-2)}=o(\varphi(\eta_{q_n}^2)).
\end{equation}
Then for any $\varepsilon>0$,
\[
P\bigl(d(g_b, \Theta_I)>\varepsilon|X^n\bigr)=o_p(1).
\]
\end{longlist}
\end{theorem}

\subsection{Normal prior}\label{sec44}
When $g_0$ is point identified, we can also establish the posterior
consistency using normal priors
%
%
\begin{equation} \label{e46}
\pi(b)=\prod_{i=1}^{q_n}\pi_i(b_i),\qquad \pi_i(b_i)\sim
N(0,\sigma^2),
\end{equation}
for some constant $\sigma^2>0$. As discussed previously, by restricting
$q_n$ to grow slowly as $n\rightarrow\infty$, we do not need a
shrinking prior to function
as a penalty term attached to the log-likelihood for the regularization
purpose.\footnote{We thank a referee for pointing this out.} Therefore~$\sigma^2$ is treated to be a fixed constant that does not depend on $n$.

With the assumptions imposed in Sections \ref{sec42} and \ref{sec43},
we can verify all
the conditions in Theorem \ref{t22}, which then leads to the following theorem:
%
%
\begin{theorem}[(Posterior consistency using Gaussian prior)] \label{t44}
Assume $g_0$ is point identified. Under Assumptions \ref{iid}, \ref
{a41}--\ref{a45}, suppose the normal prior (\ref{e46}) is used, and
%
%
\begin{equation} \label{e47}
q_n\biggl(\frac{k_n^{3d/2}}{\sqrt{n}}+\frac{1}{k_n}
\biggr)^{1/3}=o(\varphi(\eta_{q_n}^2)),
\end{equation}
then for any $\varepsilon>0$,
\[
P(\|g_b-g_0\|_s>\varepsilon|X^n)=o_p(1).
\]
\end{theorem}

\subsection{Choice of tuning parameters}

To choose $(k_n,q_n, B_n)$ that satisfy~(\ref{e44}) (\ref{e45}) and
(\ref{e47}) for each specified prior, consider the case\vadjust{\goodbreak} where $\eta
_{q_n}$ is decreasing as some power of $q_n$ [see, e.g., \citet
{Sch81} and
\citet{Mey90}], and $k_n$ grows at a polynomial rate of $n$,
that is,
%
%
\begin{eqnarray}
\eta_{q_n}&\sim&q_n^{-v}\qquad \mbox{for some }v>0,\nonumber\\[-8pt]\\[-8pt]
\frac{k_n^{3d/2}}{\sqrt{n}}+\frac{1}{k_n}&\sim& n^{-p},\qquad
0<p\leq\frac{1}{3d+2}.\nonumber
\end{eqnarray}
We then have the following corollaries:
%
%
\begin{cor}[(Truncated prior)] \label{cor41} Suppose the truncated prior
(either uniform or truncated normal) is used; then the following choice
of $(q_n, B_n)$ achieves the posterior consistency, for $b<p$:

\begin{longlist}
\item
in the mildly ill-posed case,
\[
B_n^2\sim n^b,\qquad q_n=o\bigl(n^{({p-b})/({2+2\alpha v})}\bigr);
\]
\item in the severely ill-posed case,
\[
B_n^2\sim n^b,\qquad q_n=o\bigl((\log n)^{{1}/({2\alpha v})}\bigr).
\]
\end{longlist}
\end{cor}
%
%
\begin{cor}[(Thin-tail prior)] \label{cor42} Suppose the thin-tail prior
is used; then the following choice of $q_n$ achieves the posterior
consistency, for $pr>2$:

\begin{longlist}
\item
in the mildly ill-posed case,
\[
q_n=o\bigl(n^{({pr-2})/({2r+2\alpha v(r-2)})}\bigr);
\]
\item in the severely ill-posed case,
\[
q_n=o\bigl((\log n)^{{1}/({2\alpha v})}\bigr).
\]
\end{longlist}
\end{cor}
%
%
\begin{cor}[(Normal prior)] \label{cor43} Suppose the normal prior is
used, and~$g_0$ is point identified, the following choice of $q_n$
achieves the posterior consistency:

\begin{longlist}
\item in the mildly ill-posed case,
\[
q_n=o\bigl(n^{{p}/({3(1+2\alpha v)})}\bigr);
\]

\item in the severely ill-posed case,
\[
q_n=o\bigl((\log n)^{{1}/({2\alpha v})}\bigr).
\]
\end{longlist}
\end{cor}

In the conditions of these consistency results, the
choice of tuning parameters ($q_n$, $B_n$, $r$) depend on some
parameters that one either knows or chooses ($d$, $p$), as well as
some parameters related to the true model ($\alpha$, $v$). The latter,
although undesirable, cannot be totally avoided when we study the
frequentist convergence properties under ill-posedness. [Conditions
depending on the true model are also used,
e.g., by \citet{ChePou}, directly in their Corollary 5.1, and
indirectly at the end of their Section 3.1.]

On the other hand, these results can still have meaningful implications
that do not explicitly depend on the indexes $\alpha$ and $p$ (which
are probably unknown in
practice). For example, we note that in the mildly ill-posed
situations, the
condition on $q_n$ would be satisfied if it grows as any finite power
of $\log n$. Likewise, in the severely ill-posed situations, the
condition on $q_n $
would be satisfied if it grows as any finite power of $\log\log n$.

In addition, we will indicate in the next section that the current
Bayesian-flavored treatment can even allow a data-driven choice of the
sieve dimension~$q_n$, using a posterior distribution derived from a
mixed prior.

%
%
%

\section{Random sieve dimension}\label{sec5}
As the sieve dimension $q_n$ plays an
important role not only in dealing with the ill-posed inverse problem,
but also in many applied sieve estimation methods, in this section we
briefly discuss the possibility of choosing it based on a posterior
distribution. This will require specifying a prior distribution on the
sieve dimension first. Since the conditions of a deterministic $q_n$
for consistency only restricts the growth rate, as a result, $Mq_n$
would also lead to consistency for a positive constant $M>1$, if $q_n$
ensures consistency.

We denote the sieve dimension by $q$, let it be random and place a
discrete uniform prior
%
%
\begin{equation} \label{e51}
\pi(q) = \operatorname{Unif}\{1,\ldots,Mq_n\}
\end{equation}
for some deterministic sequence $q_n\rightarrow\infty$ and constant
$M>1$. Then the prior on the sieve coefficients $b$ becomes a mixture prior
%
%
\begin{equation}\label{e52}
\pi(b)=\sum_{q=1}^{Mq_n}\pi(q)\pi(b|q)=\sum_{q=1}^{Mq_n}
(Mq_n)^{-1} \pi(b|q),
\end{equation}
where $\pi(b | q)$ follows a prior as specified before for a given
sieve dimension $q$. The feasible
limited information likelihood is, as before, denoted by $L_n(b,q)$. We
have the joint posterior
\[
p(g_b,q|X^n)\propto\pi(b|q)L_n(b,q).
\]

It can be shown that the uniform mixture prior can also lead to the
posterior consistency.
%
%
\begin{theorem}[(RANDOM $q$)] \label{t51} For each theorem in Sections
\ref{sec3} and
\ref{sec4}, suppose the corresponding conditions are satisfied for the
deterministic
sieve dimension $Mq_n$ instead of $q_n$, for some $M>1$.
Then all the posterior consistency results stated in Sections \ref{sec3}
and \ref{sec4}
(on risk consistency and on estimation consistency)
remain valid for the mixed prior (\ref{e52}) with random $q$ following
prior (\ref{e51}), with no extra conditions, with the following two
exceptions:

\begin{longlist}[(1)]
\item[(1)] We will additionally assume that $(\log q_n)/n = o(\delta_n)$ holds
for the statement of Theorem \ref{t32} to hold.

\item[(2)] We will additionally assume that $(\log q_n)/n=o(\inf_{g\in
{\mathcal
H}_n,g\notin\Theta_I^\varepsilon} G(g))$
for the statement of Theorem \ref{t32} to hold.
\end{longlist}
\end{theorem}

Note that the uniform prior is used for $q$, which gives zero prior
probability on very large choice beyond $Mq_n$. However, from a
technical point of view, the result can be extended to the case with
tails of prior on $q$ extending to infinity, as long as the tail is
thin enough so that $\pi(q>Mq_n)$ is dominated by a small enough upper bound.

The marginal posterior of $q$ is given by
%
%
\begin{equation}
p(q|X^n)\propto\int\pi(b|q)L_n(b,q)\,db.
\end{equation}
Practically, we can choose $q$ from $p(q|X^n)$.

\section{Conclusion and discussion}\label{sec6}
We studied the nonparametric conditional moment restricted model in a
quasi-Bayesian approach, with a special focus on the large sample
frequentist properties of the posterior distribution. There was no
distribution assumed on the data generating process. Instead, we
derived the posterior using the \textit{limited information likelihood}
(\textit{LIL}), allowing the proposed procedure to be simpler than the
traditional nonparametric Bayesian approach which would model the data
distribution nonparametrically. There are several alternative
moment-condition-based likelihood functions. The empirical likelihood
[\citet{Owe90})] and the generalized empirical likelihood
[\citet{ImbSpaJoh98}, \citet{NewSmi04} and \citet{Kit}]
are typical
examples. It is still possible to establish the posterior consistency
if these alternative nonparametric likelihoods are used, which is left
as a future research direction.

The parameter space $\mathcal{H}$ does not need to be compact. We
approximate $\mathcal{H}$ using a finite-dimensional sieve space
$\mathcal{H}_n$, and the regularization is carried out by a slowly
growing sieve dimension $q_n$. We then studied in detail the NPIV model
and verified all the sufficient conditions proposed in Section
\ref{sec3} in order for the posterior to be consistent.

It is also possible to achieve the posterior consistency using a larger
sieve dimension $q_n$. In this case, the regularization is carried out
by a truncated normal prior with shrinking variance, and the log-prior
is then a regularization penalty attached to the log-likelihood.
Conditions (\ref{e310}), (\ref{e311}) and Assumption \ref{a45} can
be relaxed. We describe this procedure in the Technical Report
[\citet{LiaJiaN2}].

An interesting research direction is to derive the convergence rate.
With all the tools given in this paper, it is possible to obtain the
rate of convergence of our procedure. However, the rate would be
sub-optimal, possibly due to the technical
bound (\ref{e21}) used in this paper. It would be interesting to
develop a method based on a bound tighter than (\ref{e21}), in order to
prove the nonparametric
minimax optimal rate of convergence as in \citet{ChePou09}.

In applications, our method requires a priori choices of $(k_n, q_n)$,
and $B_n$ for the truncated prior. We conjecture that the finite sample
behavior of the posterior is robust to the choice of $(k_n, B_n)$.
However, it should be sensitive to $q_n$, as a large value of $q_n$ may
lead to over-fitting. Therefore, we proposed an approach to allow for a
random sieve dimension by putting a discrete uniform prior on it and
selecting it from its posterior. With the upper bound of the uniform
prior $Mq_n$ growing under the same rate restriction as before, the
posterior consistency is also achieved. This feature, however, requires
specifying $Mq_n$. In practice, one may start with a moderate level
$Mq_n$ that is less than ten. In the NPIV setting, \citet{Hor}
recently introduced an empirical approach for selecting $q_n$.
Moreover, developing methods of selecting $(k_n, B_n)$ in a Bayesian
(or quasi-Bayesian) approach is another important research topic.



\section*{Acknowledgments}

This paper develops from a chapter of the first author's Ph.D.
dissertation at Northwestern University. We are grateful to Joel
Horowitz, Elie Tamer, Hidehiko Ichimura, Jia-Young Fu, Tom Severini,
Xiaohong Chen, Anna Simoni, an Associate Editor and two referees for
many helpful comments and suggestions on this paper. We also thank
the discussions with seminar participants at the 2010 Summer CEMMAP
conference on ``Recent developments in nonparametric instrumental
variable methods'' in London. The first author appreciates the
constant encouragements from his Ph.D. committee members at
Northwestern University.

\begin{supplement}[id=suppA]
\stitle{Technical proofs}
\slink[doi]{10.1214/11-AOS930SUPP} 
\sdatatype{.pdf}
\sfilename{aos930\_supp.pdf}
\sdescription{This supplementary material contains the proofs
of all the results developed in the main paper.}
\end{supplement}

%

\printaddresses

\end{document}